\documentclass[11pt,twosided]{article}
\usepackage{dsfont}
\usepackage{bm}
\usepackage{float}
\usepackage{srcltx}
\usepackage{mathrsfs}
\usepackage{graphicx}
\usepackage{amsfonts}
\usepackage{indentfirst}
\usepackage{amssymb,amsfonts,amsmath, 
colordvi}
\usepackage{diagbox}
\usepackage{caption}

\usepackage{leftidx}
\usepackage{pgf,tikz}
\usepackage{cite}
\usetikzlibrary{shapes,arrows}
\allowdisplaybreaks[4]
\usepackage{color}

\usepackage{amssymb,color}

\definecolor{c20}{rgb}{0.,0.7,0.}
\definecolor{c30}{rgb}{0.,0.,1.}
\definecolor{c40}{rgb}{1,0.1,0.7}
\definecolor{c50}{rgb}{1,0,0}
\definecolor{c60}{rgb}{0,0.9,0.1}

\usepackage{amssymb,color}

\newcommand{\N}{\mathbb{N}}

\newcommand{\BQN}{\begin{eqnarray}}
\newcommand{\EQN}{\end{eqnarray}}
\newcommand{\BQNY}{\begin{eqnarray*}}
\newcommand{\EQNY}{\end{eqnarray*}}

\newcommand{\BS}{\begin{sat}}
\newcommand{\ES}{\end{sat}}
\newcommand{\BT}{\begin{theo}}
\newcommand{\ET}{\end{theo}}
\newcommand{\BK}{\begin{korr}}
\newcommand{\EK}{\end{korr}}

\newcommand{\BD}{\begin{de}}
\newcommand{\ED}{\end{de}}

\newcommand{\BIT}{\begin{itemize}}
\newcommand{\EIT}{\end{itemize}}
\newcommand{\BDI}{\begin{description}}
\newcommand{\EDI}{\end{description}}

\newcommand{\BRM}{\begin{remarks}}
\newcommand{\ERM}{\end{remarks}}

\newcommand{\BEL}{\begin{lem}}
\newcommand{\EEL}{\end{lem}}

\newtheorem{theo}{Theorem}[section]

\newtheorem{lem}[theo]{Lemma}

\newtheorem{coro}[theo]{Corollary}
\newtheorem{re}[theo]{Remark}
\newtheorem{remarks}[theo]{Remarks}

\newcommand{\COM}[1]{}

\begin{document}
\title{Optimal Parisian-type dividends payments discounted by the number of claims for the perturbed classical risk process}
\author{Irmina Czarna\thanks{Mathematical Institute, University of Wroclaw, Poland. E-mail: czarna@math.uni.wroc.pl}\ \ Yanhong Li\thanks{School of mathematical Sciences and LPMC Nankai University, Tianjin 300071, P.R. China. E-mail: yanhongli@mail.nankai.edu.cn}\ \ \ Zbigniew Palmowski\thanks{Mathematical Institute, University of Wroclaw, Poland. E-mail: zbigniew.palmowski@gmail.com}  \ \ Chunming Zhao\thanks{Department of Statistics
School of Mathematics Southwest Jiaotong University, Chengdu, Sichuan, 611756, PR of China. Email:cmzhao@swjtu.cn}}

\date{}
\maketitle
\begin{abstract}
In this paper we consider a classical risk process perturbed by a Brownian motion.
We analyze the value function describing the mean of the cumulative discounted dividend payments paid up to Parisian ruin time and further discounted by the number of claims appeared up to this ruin time. We identify this value function for the
barrier strategy and find the sufficient conditions for this strategy to be optimal.
We also consider few particular examples.
\end{abstract}


\section{Introduction}
The classical optimal dividends problem has been been considered by many authors since De Finetti \cite{de1957impostazione} who introduced it
to address the objection that risk process has the unrealistic property that it converges to
infinity with probability one.
Gerber and Shiu \cite{gerber2004optimal}, Asmussen and Taksar \cite{asmussen1997controlled} and Jeanblanc and Shiryaev \cite{jeanblanc1995optimization} considered the the optimal dividend problem in the Brownian motion setting. Azcue and Muler \cite{azcue2005optimal} and Schmidli \cite{schmidli2006optimisation} studied the optimal dividend strategy under the Cram\'{e}r-Lundberg model using a Hamilton-Jacobi-Bellman (HJB) system of equations. Further Avram et.al. \cite{avram2007optimal, avram2015gerber}, Kyprianou and Palmowski \cite {kyprianou2007distributional}, Loeffen \cite{loeffen2008optimality, loeffen2009optimal}, Loeffen and Renaud \cite{loeffen2010finetti}, Czarna and Palmowski \cite{czarna2010dividend} and many other authors analyze the divident problem
for the L\'{e}vy risk process using the probabilistic approach.

In this paper will consider a classical risk process perturbed by a Brownian motion, which a particular case of L\'evy risk process
covering vast range of examples which are of insurance companies interest. The classical risk process models so-called 'large' claim size
and independent Brownian motion models 'small' claims.

In this paper we want to combine two new ingredients appearing in the value function.
The first one concerns choosing the Parisian ruin time instead of classical ruin time.
Parisian ruin occurs if the process falls below zero and stays below zero for a continuous time interval of length $d\geq $.
The name comes from a Parisian option which is activated or canceled depending on if the underlying asset price stays above or below the barrier for a long enough period of time or not (see \cite{chesney1997brownian} and \cite{dassios2008parisian}).
We believe that giving a possibility of Parisian delay could be better in many situations as it gives the insurance companies the chance to achieve solvency.
Still, the particular case of $d=0$ brings us to the classical set-up.
The second new factor in the objective function concerns additional discounting component related with
the total number of claims arrived until Parisian ruin.
Taking into account extra discounting based on the total number of the claims/losts arrived up to ruin time allows to diminish the objective value in the case of
large number of claims and increase it in the case of the small number.
This is a very natural and practical feature of dividends payments.
Therefore we believe that all analysis presented in this paper will be vary
valuable for the actuarial applications.

We identify above described value function for the
barrier strategy according which any surplus above fixed level $a$ is paid as dividends.
We also find the sufficient conditions for this strategy to be optimal.
Finally, we calculate the value function for few examples.

The paper is organized as follows. In Sections 2 we introduce the basic notions, notations and describe in detail the model we deal in this paper with.
In Section 3 we analyze the barrier strategy. In Section 4 we find the sufficient conditions under which the barrier strategy is optimal.
In last section we present the example when the claim size is exponential distributed.

\section{Risk process and value function}
We consider the following surplus process in continuous time:
\begin{equation}\label{Xt}
   X_t=x+ct-\sum_{i=1}^{N_t}C_{i}+\sigma B_t,
\end{equation}
where the non-negative constant $x$ denotes the initial reserve (later we underline this initial capital by adding subscript to probability measures $\mathbb{P}_x$ with $\mathbb{P}:=\mathbb{P}_0$ and to
corresponding expectations $\mathbb{E}_x$  with $\mathbb{E}:=\mathbb{E}_0$).
$N_t$ is a Poisson process with intensity $\lambda$ describing the number of claims appeared
up to time $t$, $\{C_{i}\}_{i=1}^{\infty}$ are claim sizes which are independent and identically distributed  (i.i.d.) non-negative  random variables that are also independent of $N_t$. By $f(x)$ and $F(x)$ we denote the density and distribution function of generic $C$, respectively.
Throughout of this paper we assume that the claim density $f$ is continuously differentiable.
The positive constant $c=\lambda\mathbb{E}(C_{1})(1+\theta)$ is the rate of premium income and $\theta>0$ is the relative security loading factor.
$B_t$ is a standard Wiener process that is independent of the aggregate claims process $\sum_{i=1}^{N_t}C_{i}$ and $\sigma\geq 0$ is a dispersion parameter. See also Figure 1 for an example of sample path where $\sigma=0$.
\begin{center}
\definecolor{ffqqtt}{rgb}{1,0,0.2}
\begin{tikzpicture}[line cap=round,line join=round,>=triangle 45,x=1.0cm,y=1.0cm]
\draw[->,color=black] (-1,0) -- (11.5,0);
\foreach \x in {}
\draw[shift={(\x,0)},color=black] (0pt,-2pt);
\draw[->,color=black] (0,-2) -- (0,4);
\foreach \y in {}
\draw[shift={(0,\y)},color=black] (2pt,0pt) -- (-2pt,0pt);
\clip(-1,-3) rectangle (12,4);
\draw[dash pattern=on 3pt off 3pt,color=ffqqtt](0,1.5)-- (11.5,1.5);
\draw (0,1)-- (1,2);
\draw [dash pattern=on 3pt off 3pt] (1,2)-- (1,0.5);
\draw (1,0.5)-- (2,1.5);
\draw [dash pattern=on 3pt off 3pt] (2,1.5)-- (2,-1);
\draw (2,-1)-- (5,2);
\draw [dash pattern=on 3pt off 3pt] (5,2)-- (5,1);
\draw (5,1)-- (7,3);
\draw [dash pattern=on 3pt off 3pt] (7,3)-- (7,-2);
\draw (7,-2)-- (8,-1);
\draw [dash pattern=on 3pt off 3pt] (8,-1)-- (8,-1.5);
\draw (8,-1.5)-- (9,-0.5);
\draw [dash pattern=on 3pt off 3pt] (9,-0.5)-- (9,-1);
\draw (9,-1)-- (10.5,0.5);
\draw [shift={(3,-3.26)},color=ffqqtt]  plot[domain=1.27:1.87,variable=\t]({1*3.41*cos(\t r)+0*3.41*sin(\t r)},{0*3.41*cos(\t r)+1*3.41*sin(\t r)});
\draw [shift={(8,-3.26)},color=ffqqtt]  plot[domain=1.27:1.87,variable=\t]({1*3.41*cos(\t r)+0*3.41*sin(\t r)},{0*3.41*cos(\t r)+1*3.41*sin(\t r)});
\draw (-0.95,4) node[anchor=north west] {$X_t$};
\draw (11,0.07) node[anchor=north west] {$t$};
\draw (-0.3,0.1) node[anchor=north west] {$0$};
\draw (-0.53,1.30) node[anchor=north west] {$u$};
\draw (-0.53,1.70) node[anchor=north west] {$a$};
\draw (2.82,0.79) node[anchor=north west] {$d$};
\draw (8.04,0.84) node[anchor=north west] {$d$};
\draw (0.5,-2.5) node[anchor=north west] {{\rm Figure\; 1.\; A sample path of the original  surplus process\; $X_t.$}};
\end{tikzpicture}
\end{center}

We denote a dividend or control strategy by $\pi$, where $\pi=\{L^{\pi}_t: t\geq0\}$ is a nondecreasing left-continuous adapted process which starts at zero. The random variable $L^{\pi}_t$ represents the cumulative dividends the company has paid out up to time $t$ under the control $\pi$. We define by
\begin{equation}\label{Upi}
U^{\pi}_t=X_t-L^{\pi}_t
\end{equation}
the controlled risk process under the dividend strategy $\pi$. For fixed $d\geq 0$ let
\begin{equation}\label{ruintime}
 \tau^{\pi,d}:=\inf\{t>0: t-\sup\{s<t: U^\pi_s\ge0\}> d, U^\pi_t<0\}
\end{equation}
be a Parisian ruin time, that is, the ruin occurs if the regulated process $U^\pi_t$ falls below zero
and stays below zero for a continuous time interval of length $d$ (see Figure 1).
The case $d=0$  corresponds to the classical ruin time:
$$\tau^\pi:=\inf\{t\geq 0: U^\pi_t<0\}.$$
Now we formally define the value function for a dividend strategy $\pi$ as follows:
\begin{equation*}\label{v}
 v^{\pi}(x):=\mathbb{E}_x[r^{N_{\tau^{\pi,d}}}\int^{\tau^{\pi,d}}_0e^{-qt}dL^{\pi}_t],
\end{equation*}
where $r\in(0,1]$ is a constant and $q>0$ is a discounting rate.
By definition, it follows that
\begin{equation}\label{zero}
v^{\pi}(x)=0\qquad \text{for $x<-cd$}
\end{equation}
since the risk process will not manage to come back to zero before Parisian ruin and we will
collect no dividends. A strategy $\pi$ is called admissible if ruin does not occur by a dividend payout, that is, $L^{\pi}_{t+}-L^{\pi}_t\leq U^{\pi}_t$ for $t<\tau^{\pi}$. Let $\Pi$ denote the set of all admissible dividend strategies. The objective of the beneficiaries of the insurance company is to maximize $v^{\pi}(x)$ over all admissible strategies $\pi$, that is, finding the optimal value function $v^*$ given by
\begin{equation*}
v^{*}(x)=\sup_{\pi\in\Pi}v^{\pi}(x)
\end{equation*}
and identifying the optimal strategy $\pi^{*}\in\Pi$ such that
 \begin{equation*}
v^{\pi^*}(x)=v^{*}(x)        \ \ \ \ \ \ \ \ \ \ \ \ \ \ \ \ \ \text{for  all $x> -cd$.}
\end{equation*}

The crucial dividends distribution policy  is the {\it barrier
policy} $\pi_a=\{L^a_t:=L^{\pi_a}_t: t\geq0\}$ of transferring all surpluses above a given level
$a$ as dividends; see Figure 2. In this case $L_t^a:=a\wedge \sup_{s\leq t}X_s-a$.
For this specific strategy we will denote
\begin{equation}\label{Ut}
U^a_t:=U^{\pi_a}_t=X_t-L^a_t,\qquad \tau^d:=\tau^{\pi_a,d}
\end{equation}
with
\begin{equation}\label{taud}\tau:=\tau^{0}\end{equation}
being classical ruin time
and
\begin{equation}\label{v}
v^{a,d}(x):=v^{\pi_a,d}=\mathbb{E}_x[r^{N_{\tau^d}}\int^{\tau^d}_0e^{-qt}dL^a_t].
\end{equation}
In this paper we will find the value function $v^{a,d}(x)$ under the barrier strategy $\pi_a$.
Later we will identify the optimal barrier $a^*$ maximizing
$v^{a,d}(x)$ for any $x>-cd $. Finally in the last step we show that
in most of the known cases of the claim size density $f$
the strategy $\pi_{a^*}$ is indeed optimal, that is
\begin{equation}\label{optimality}
v^*(x)=v^{a^*,d}(x),\qquad x>-cd.
\end{equation}

\begin{center}
\definecolor{ffqqtt}{rgb}{1,0,0.2}
\begin{tikzpicture}[line cap=round,line join=round,>=triangle 45,x=1.0cm,y=1.0cm]
\draw[->,color=black] (-1,0) -- (11.5,0);
\foreach \x in {}
\draw[shift={(\x,0)},color=black] (0pt,-2pt);
\draw[->,color=black] (0,-3.5) -- (0,3);
\foreach \y in {}
\draw[shift={(0,\y)},color=black] (2pt,0pt) -- (-2pt,0pt);
\clip(-1,-3,5) rectangle (12,3);
\draw[dash pattern=on 3pt off 3pt,color=ffqqtt](0,1.5)-- (11.5,1.5);
\draw (0,1)-- (0.5,1.5);
\draw (0.5,1.5)-- (1,1.5);
\draw [dash pattern=on 3pt off 3pt] (1,1.5)-- (1,0.1);
\draw (1,0.1)-- (2,1.1);
\draw [dash pattern=on 3pt off 3pt] (2,1.1)-- (2,-1.4);
\draw (2,-1.4)-- (4.9,1.5);
\draw (4.9,1.5)-- (5,1.5);
\draw [dash pattern=on 3pt off 3pt] (5,1.5)-- (5,0.5);
\draw (5,0.5)-- (6,1.5);
\draw (6,1.5)-- (7,1.5);
\draw [dash pattern=on 3pt off 3pt] (7,1.5)-- (7,-3.5);
\draw (7,-3.5)-- (8,-2.5);
\draw [dash pattern=on 3pt off 3pt] (8,-2.5)-- (8,-3);
\draw (8,-3)-- (9,-2);
\draw [dash pattern=on 3pt off 3pt] (9,-2)-- (9,-2.5);
\draw (9,-2.5)-- (10.5,-1);
\draw [shift={(3,-3.26)},color=ffqqtt]  plot[domain=1.27:1.87,variable=\t]({1*3.41*cos(\t r)+0*3.41*sin(\t r)},{0*3.41*cos(\t r)+1*3.41*sin(\t r)});
\draw [shift={(8,-3.26)},color=ffqqtt]  plot[domain=1.27:1.87,variable=\t]({1*3.41*cos(\t r)+0*3.41*sin(\t r)},{0*3.41*cos(\t r)+1*3.41*sin(\t r)});
\draw (-0.95,3) node[anchor=north west] {$U^a_t$};
\draw (11,0.07) node[anchor=north west] {$t$};
\draw (-0.56,0.14) node[anchor=north west] {$0$};
\draw (-0.53,1.30) node[anchor=north west] {$u$};
\draw (-0.53,1.70) node[anchor=north west] {$a$};
\draw (2.82,0.79) node[anchor=north west] {$d$};
\draw (8.04,0.84) node[anchor=north west] {$d$};
\draw (8.86,0.15) node[anchor=north west] {$\tau^d$};
\draw (1.9,0) node[anchor=north west] {$\tau$};
\draw (0.5,-4) node[anchor=north west] {${\rm Figure \; 2.\;  A \; sample\; path\; of\; the\; regulated\; surplus\; process}\; U^a_t.$};
\end{tikzpicture}
\end{center}

\section{Barrier strategy $\pi_a$}

\subsection{The value function}

We start from identifying the value function $v^{a,d}(x)$ under the barrier strategy $\pi_a$ for a fixed $a$.
We will prove now basic properties of $v^{a,d}(x)$.

\begin{lem}\label{lemeq1}
Function $v^{a,d}$ is continuously differentiable for all $x<a$ and solves the following differential equation:
\begin{equation}\label{hidea}
\frac{\sigma^2}{2}(v^{a,d})^{''}(x)+c(v^{a,d})^{'}(x)-(\lambda+q)v^{a,d}(x)+\lambda r\int_0^{\infty}v^{a,d}(x-y)f(y)dy=0.
\end{equation}
\end{lem}

\textbf{Proof.}
The proof is based on the Theorem 3 of Li \cite{li2006distribution}.
Therefore we decided to skip all details.
Still it requires really few modifications.
The main one is in the statement of Theorem 2 of  of Li \cite{li2006distribution}, where
we should write $u+cs+\sigma y+cd$ or simply to $\infty$ (remembering (\ref{zero}))
instead of $u+cs+\sigma y$ in the integral of the second increment.
Apart of this, $\lambda$ should be accompanied by additional $r$.
The difference comes from the fact that Parisian ruin allows the regulated risk process goes below zero but
not below $-cd$ since otherwise there will be ruin after delay $d$.
Moreover, when the first claim arrives, the discounting with respect of number of arrived claims gives additional $r$.
This small change should be continued in all calculations regarding $I_3(u)$.
This completes the proof.

\vspace{3mm}\hfill $\Box$

We will express the value function $v^{a,d}(x)$ in terms of some special function
$h^d(x)$ related with the first passage times which we will now formally define.
For $y\in\mathbb{R}$, we define the first passage times of $X$:
$$
\tau^{+}_y: =\inf\{t\geq0: X_t\geq y\}, \qquad \tau^{-}_y: =\inf\{t\geq0: X_t< y\}.
$$

For $-cd\leq x\leq a$ we define
 \begin{equation}\label{hddefinition}
h^d(x):=\mathbb{E}_x[r^{N_{\tau^+_a}}e^{-q\tau^+_a}, \tau^+_a<\tau^d]
\end{equation}
and for $x<-cd$
\begin{equation*}
h^d(x)=0.
\end{equation*}

Observe that our risk process given in (\ref{Xt}) is a spectrally negative L\'evy process.
Moreover, it is a Markov process.
Hence, using fact that this process up-crosses all levels continuously, from the Markov property it follows that
for general $b>a$
$$\mathbb{E}_x[r^{N_{\tau^+_b}}e^{-q\tau^+_b}, \tau^+_b<\tau^d]=
\mathbb{E}_x[r^{N_{\tau^+_a}}e^{-q\tau^+_a}, \tau^+_a<\tau^d]\mathbb{E}_a[r^{N_{\tau^+_b}}e^{-q\tau^+_b}, \tau^+_b<\tau^d]$$
Taking $g^d(x):=\mathbb{E}_x[r^{N_{\tau^+_b}}e^{-q\tau^+_b}, \tau^+_b<\tau^d]$
we have:
\begin{equation}\label{twosidedexit}
h^d(x)=\frac{g^d(x)}{g^d(a)}.
\end{equation}
In fact for any $x<y$ we have:
\begin{equation}\label{twosidedexitb}
\mathbb{E}_x[r^{N_{\tau^+_y}}e^{-q\tau^+_y}, \tau^+_y<\tau^d]=\frac{g^d(x)}{g^d(y)}.
\end{equation}

For $x\leq a$ from the Markov property applied at the stopping time $\tau_a^+$ to $X_t$ and for $x>a$
straightforward from the definition of the barrier strategy we can conclude the following fact.
\begin{lem}\label{lemma1}
The value function $v^{a,d}(x)$ calculated for the barrier strategy $\pi_a$ satisfies the following equations:
 \begin{equation}\label{xleqad}
v^{a,d}(x)=h^d(x)v^{a,d}(a)\qquad \text{for $x\leq a$}
\end{equation}
and
\begin{equation}\label{xgeqad}
v^{a,d}(x)=x-a+v^{a,d}(a)  \qquad \text{for $x>a$.}
\end{equation}
\end{lem}

\begin{re}\rm
Note that by Lemmas \ref{lemeq1} and \ref{lemma1} function $h^d$ similarly like $v^{a,d}(x)$ is twicely differantiable and
it solves for $x\leq a$ the same differential equation:
\begin{equation}\label{hide}
\frac{\sigma^2}{2}(h^{d})^{''}(x)+c(h^{d})^{'}(x)-(\lambda+q)h^{d}(x)+\lambda r\int_0^{\infty} h^{d}(x-y)f(y)dy=0
\end{equation}
with the obvious boundary condition
\begin{equation}\label{bhd}
h^d(a)=1.\end{equation}
\end{re}

We can now conclude the following representation of the value function.

\begin{theo}\label{jieguoteo}
The value function of the barrier strategy at level $a\geq0$ is given by:
\begin{equation}\label{jieguo}
v^{a,d}(x)=\left\{
  \begin{array}{lll}
&\frac{h^d(x)}{(h^d)^{'}(a)}&\text{if}\ x\leq a,\\
&x-a+\frac{1}{(h^d)^{'}(a)}&\text{if} \ x> a.
  \end{array}
\right.
\end{equation}
\end{theo}

\begin{re}\rm
Note that the function $v^{a,d}(x)$ is differentiable for all $x\geq 0$ and
\begin{equation}\label{vad1}
(v^{a,d})^{'}(a)=1.
\end{equation}
\end{re}

\textbf{Proof of Theorem \ref{jieguoteo}.}
Note that for $n\in \N$,
\begin{eqnarray*}v^{a,d}(a)&\geq& \mathbb{E}_a\left[r^{N_{\tau_{a+1/n}^+}}e^{-q\tau_{a+1/n}^+}, \tau_{a+1/n}^+<\tau^d\right]v^{a,d}\left(a+\frac{1}{n}\right)\\&=&\mathbb{E}_a\left[r^{N_{\tau_{a+1/n}^+}}e^{-q\tau_{a+1/n}^+}, \tau_{a+1/n}^+<\tau^d\right]\left(v^{a,d}(a)+\frac{1}{n}\right)\end{eqnarray*}
and
\begin{eqnarray*}v^{a,d}(a)&\leq& \mathbb{E}_a\left[r^{N_{\tau_{a+1/n}^+}}e^{-q\tau_{a+1/n}^+}, \tau_{a+1/n}^+<\tau^d\right]\left(v^{a,d}(a)+\frac{1}{n}\right)\\&&
+\frac{1}{n}\mathbb{E}_a\left[r^{N_{\tau_{a+1/n}^+}}\int_0^{\tau_{a+1/n}^+}e^{-qt}dt, \tau_{a+1/n}^+<\tau^d\right]\\
&&+\mathbb{E}_a\left[r^{N_{\tau^d}}\int_0^{\tau^d}e^{-qt}dL^a_t, \tau^d<\tau_{a+1/n}^+\right]
\end{eqnarray*}
since $L_t^{a}=\overline{X}_t-a$ under $\mathbb{P}_a$ can increase only by $1/n$ each time the regulated process is above $a$
up to time $\tau_{a+1/n}^+$.
Moreover, since $r\leq 1$ we have
\[\mathbb{E}_a\left[r^{N_{\tau^d}}\int_0^{\tau^d}e^{-qt}dL^a_t, \tau^d<\tau_{a+1/n}^+\right]\leq \frac{1}{n}\mathbb{P}_a\left(\tau^d<\tau_{a+1/n}^+\right).\]
Note that
$\lim_{n\to\infty}\mathbb{E}_a\left[r^{N_{\tau_{a+1/n}^+}}\int_0^{\tau_{a+1/n}^+}e^{-qt}dt, \tau_{a+1/n}^+<\tau^d\right]=0$
and $\lim_{n\to\infty}\mathbb{P}_a\left(\tau^d<\tau_{a+1/n}^+\right)=0$.
All details of above estimation could be also found in the proof of Proposition 1 of Renauld and Zhou \cite{JFZ} where
$T_n^\prime=\min\{\tau_{a+1/n}^+, \tau^0\}$ should be exchanged into $T_n^\prime=\min\{\tau_{a+1/n}^+, \tau_0^d\}$.
Taking in definition of $g_d(x)$ parameter $b>a+1/n$  from (\ref{twosidedexit}) we will derive that
$$v^{a,d}(a)= \frac{g^d(a)}{g^d\left(a+\frac{1}{n}\right)}\left(v^{a,d}(a)+\frac{1}{n}\right)+{\rm o}\left(\frac{1}{n}\right).$$
Hence
\begin{equation}\label{vafinal}v^{a,d}(a)=\frac{g^d(a)}{g^{d,\prime}(a)}=\frac{1}{(h^{d})^{'}(a)}.\end{equation}
The assertion of the theorem follows now from equations (\ref{xleqad}) and  (\ref{vafinal}).
\vspace{3mm}\hfill $\Box$

Obviously to maximize the value function under barrier strategy we should choose level $a$ that minimizes $(h^d)^{'}(a)$.
\begin{theo}
The optimal barrier $a^*$ such that $v^{a^*,d}(x)=\max_{a\geq 0}v^{a,d}(x)$ solves:
\begin{equation*}
(h^d)^{''}(a^*)=0.
\end{equation*}
\end{theo}

\subsection{Identification of function $h^d(x)$}

We will now identify function $h^d(x)$. This allows to find the optimal value function and optimal strategy.

Let
\begin{equation*}
 v_y(k,t)=\frac{d}{dt}V_y(k,t),\ \ \ k\in\mathbb{N},\quad t\ge y/c
 \end{equation*}
  with
\begin{equation*}
V_y(k,t)=\mathbb{P}(N_{\tau_y^+}=k, \tau_y^+\leq t), \ \ \ k\in\mathbb{N},\quad t\ge 0.
 \end{equation*}
For $r\in(0,1]$ and $q>0$  we define the joint Laplace transform of $N_{\tau_y^+}$ and $\tau_y^+$:
\begin{eqnarray}\label{define phi}
    \phi(y)&:=&\mathbb{E}[r^{N_{\tau_y^+}}e^{-q \tau_y^+}\mathbb{I}(\tau_y^+<\infty)]=\int_0^\infty e^{-q t}\sum_{k=0}^\infty r^k v_y(k,t)dt.
\end{eqnarray}

\begin{lem}\label{v y k tlem}
If $\sigma >0$ then
\begin{eqnarray}\label{v y k t}
v_y(k,t)&=&\frac{\lambda^k}{k!}yt^{k-1}e^{-\lambda t}\int_{-\infty}^{\infty}\frac{1}{\sqrt{2\pi t}}e^{-\frac{x^2}{2t}}f^{k*}(ct+\sigma x-y)dx,
\end{eqnarray}
where $f^{k*}$ denotes $k$th convolution of the claim density function.
If $\sigma=0$ then
 \begin{eqnarray}\label{v y k t0}
v_y(k,t)&=&\frac{\lambda^k}{k!}yt^{k-1}e^{-\lambda t}f^{k*}(ct-y).
\end{eqnarray}
\end{lem}
{\bf Proof.} We assume firstly that $\sigma >0$.
Similarly like in the case of proving Lemma \ref{lemeq1} one can show that
$\phi$ is twicely-differentiable and it solves the following integro-differential equation:
\begin{equation}\label{integro differential equation}
\frac{\sigma^2}{2} \phi^{''}(y)-c\phi^{'}(y)-(\lambda+\delta)\phi(y)+\lambda r\int_0^{\infty}\phi(y+x)f(x)dx=0.
\end{equation}
Clearly, we have two boundary conditions:
\begin{equation}\label{boundary condition}
\phi(0)=1
\end{equation}
and
\begin{equation}\label{boundary condition 2}
\lim_{y\to\infty}\phi(y)=0.
\end{equation}
Since the solution to (\ref{integro differential equation}) with the boundary conditions (\ref{boundary condition}) and (\ref{boundary condition 2}) is unique,
we will check that $\phi(y)$ is of the form
\begin{equation*}
\phi(y)=e^{-b_0y}
\end{equation*}
for some $b_0$.
Note that the real part of $b_0$ must be positive, because otherwise it would be a contradiction to the fact that $\lim_{y\rightarrow\infty}\phi(y)=0$.
Let $\rho$ be a unique solution to the fundamental Lundberg's equation (see Gerber and Shiu (1997)\cite{gerber1997joint}):
\begin{equation}\label{lundberg}
 \frac{\sigma^2}{2}s^2+cs-(\lambda+q)+\lambda r\hat{f}(s)=0.
\end{equation}
Observe that:
\begin{equation}
\phi(y)=e^{-\rho y}.
\end{equation}
Taking inverse Laplace transform with respect of $q$ it follows from (\ref{define phi}) that
\begin{eqnarray}\label{e rho y}
\sum_{k=0}^{\infty}r^kv_y(k,t)&=&\mathcal{L}^{-1}_{q}(e^{-\rho y})\nonumber\\
&=&\sum_{k=0}^{\infty}r^k \frac{\lambda^k}{k!}yt^{k-1}e^{-\lambda t}\int_{-\infty}^{\infty}\frac{1}{\sqrt{2\pi t}}e^{-\frac{x^2}{2t}}f^{k*}(ct+\sigma x-y)dx.
\end{eqnarray}
Details can be found in Czarna et al. \cite{Czarna2015number}.
The case of $\sigma=0$ could be proved in a similar way.
\vspace{3mm}\hfill $\Box$

\begin{re}\rm
In particular, when the claim amounts are exponentially distributed with mean $1/\mu$ then
\begin{eqnarray}\label{expv}
v_y(k,t)=\left\{
  \begin{array}{lll}
&\frac{y}{t}e^{-\lambda t}\delta_0(ct-y),&k=0,\\
&\frac{\lambda^k \mu^{k}y}{k!(k-1)!}t^{k-1}(ct-y)^{k-1}e^{-(\lambda+\mu c)t}  e^{\mu y},&k>0.
  \end{array}
\right.
\end{eqnarray}
\end{re}

\bigskip

The function $h^d(x)$ will be given using the Dickson operator $T_r$ for $ \mathfrak{Re}(r)\geq0$ (see Dickson and Hipp \cite{dickson2001time}).
For any integrable real-valued function $g$ it is defined as
\begin{equation*}
T_r g(x):=\int_x^{\infty}{e^{-r(u-x)}g(u)}du,\ \   x\geq0.
\end{equation*}
The operator $T_{r}$ satisfies the following properties:
\begin{enumerate}
  \item $ \hat{g}:=T_s g(0)=\int_0^{\infty}{e^{-s x}g(x)dx}=\hat g(s)$ which is the Laplace transform of $g$;
\item  $T_rT_s g(x)=\int_x^{\infty}e^{-r(y-x)}\int_y^{\infty}e^{-s(z-y)}g(z)dzdy$.
\item The operator $T_s$ is commutative, i.e. $T_sT_r=T_rT_s$. Moreover,
\begin{equation}\label{Dickson's Operator}
T_{s}T_{r}f(x)=T_{r}T_{s}f(x)=\frac{T_s f(x)-T_r f(x)}{r-s}, s\neq r, x\geq0.
\end{equation}
\end{enumerate}
More properties of the operator $T_r$ can be found in Li and Garrido \cite{li2004ruin}.

Denote:
\begin{equation}\label{wd}
w_d(x):=\int_0^{\infty}\int_{0}^de^{-qt}\sum_{k=0}^{\infty}r^kv_z(k,t)dtf(y+x)dy,
\end{equation}
which by Lemma \ref{v y k tlem} could be easily identified.
Next theorem gives explicit expression for the function $h_d(x)$ producing the value function $v^{a,d}(x)$.

\begin{theo}
For $0\leq x\leq a$ the function $h^d(x)$ can be expressed as follows:
\begin{equation}\label{hxr}
h^d(x)=\frac{\sum_{n=0}^{\infty}(\frac{\lambda r}{c})^n(T_{\rho}f)^{*n}*\zeta(x)}{\sum_{n=0}^{\infty}(\frac{\lambda r}{c})^n(T_{\rho}f)^{*n}*\zeta(a)},\ \ \ \ \  \ \ \ \ \ \ \ \text{for $\sigma=0,$,}\\
\end{equation}
 and
 \begin{equation}\label{hdxr}
 h^d(x)=\frac{\sum_{n=0}^{\infty}(\frac{2\lambda r}{\sigma^2})^{n}(\beta*T_{\rho}p)^{*n}*\varphi(x)}{\sum_{n=0}^{\infty}(\frac{2\lambda r}{\sigma^2})^{n}(\beta*T_{\rho}p)^{*n}*\varphi(a)},\ \ \ \ \text{for $\sigma>0,$,}
 \end{equation}
 where $$\varphi(x):=[\rho+2c/\sigma^2+\xi^{'}(0)]f*\beta(x)+\beta(x)-\frac{2\lambda r}{\sigma^2}f*\beta*w_d(x),$$
  $$\beta(x):=e^{-(\rho+\frac{2c}{\sigma^2})x},\qquad \zeta(x):=e^{\rho x}.$$
 \end{theo}

{\bf Proof.}
Let at the beginning $\sigma >0$ and $d>0$.
Take any $y<0$. Then form the Markov property we have:
\begin{eqnarray*}
h^d(y)&=&\mathbb{E}_y[r^{N_{\tau_0^{+}}}e^{-q\tau_0^{+}}, \tau_0^{+}<d]h^d(0)\nonumber\\
&=&\mathbb{E}_0[r^{N_{\tau_{-y}^{+}}}e^{-q\tau_{-y}^{+}}, \tau_{-y}^{+}<d]h^d(0)\nonumber\\
&=&h^d(0)\int_0^de^{-qt}\sum_{k=0}^{\infty}r^kv_{-y}(k,t)dt,
\end{eqnarray*}
Then the equation (\ref{hide}) is equivalent to the following equation:
\begin{eqnarray}\label{hideb}
\lefteqn{\frac{\sigma^2}{2}(h^{d})^{''}(x)+c(h^{d})^{'}(x)-(\lambda+q)h^d(x)+\lambda r\int_0^{x}h^d(x-y)f(y)dy}\nonumber\\&&+\lambda rh^d(0)\int_x^{\infty}f(y)\int_0^de^{-qt}\sum_{k=0}^{\infty}r^kv_{y-x}(k,t)dtdy=0.
\end{eqnarray}
Since $\mathbb{P}(\tau^{+}_a<\tau^d)>0$ the definition of $h^d(x)$ yields that $h^d(0)\ne0$.
Dividing both sides in (\ref{hideb}) by $h^d(0)$ and letting $\xi(x)=\frac{h^d(x)}{h^d(0)}$ we produce the following equation:
\begin{equation}\label{gb}
\frac{\sigma^2}{2}\xi^{''}(x)+c\xi^{'}(x)-(\lambda+q)\xi(x)+\lambda r\int_0^{x}\xi(x-y)f(y)dy+\lambda rw_d(x)=0,
\end{equation}
where $w_d(x)$ is given in (\ref{wd}).
Taking Laplace transforms on both sides of (\ref{gb}) for sufficiently large $s$ we get
   \begin{equation}\label{gl}
   \frac{\sigma^2}{2}s^2\hat{\xi}(s)-\frac{\sigma^2}{2}\xi(0)s-\frac{\sigma^2}{2}\xi^{'}(0)+cs\hat{\xi}(s)-c\xi(0)-(\lambda+q)\hat{\xi}(s)+\lambda r\hat{f}(s)\hat{\xi}(s)+\lambda r\hat{w_d}(s) =0
   \end{equation}
which after rearranging terms leads to
  \begin{equation}\label{gla}
    \left[ \frac{\sigma^2}{2}s^2+cs-(\lambda+q)+\lambda r\hat{f}(s)\right]\hat{\xi}(s)=\frac{\sigma^2}{2}\xi(0)s+\frac{\sigma^2}{2}\xi^{'}(0)+c\xi(0)-\lambda r\hat{w_d}(s) .
   \end{equation}
 By the definition of $\xi(x)$ it follows that $\xi(0)=1$. Then (\ref{gla}) gives:
\begin{equation}\label{gla1}
    \left[\frac{\sigma^2}{2}s^2+cs-(\lambda+q)+\lambda r\hat{f}(s)\right]\hat{\xi}(s)=\frac{\sigma^2}{2}s+\frac{\sigma^2}{2}\xi^{'}(0)+c-\lambda r\hat{w_d}(s) .
   \end{equation}
Recall that $\rho$ is the unique nonnegative root of the Lundberg fundamental equation (\ref{lundberg}). Then
subtracting $[\frac{\sigma^2}{2}\rho^2+c\rho-(\lambda+q)+\lambda r\hat{f}(\rho)]\hat{\xi}(s)$ on the left side of (\ref{gla1})
and dividing the result equation by $s-\rho$ produces
\begin{equation}\label{glap}
     \left[\frac{\sigma^2}{2}(s+\rho)+c-\lambda rT_sT_\rho f(0)\right]\hat{\xi}(s)=\frac{\sigma^2}{2}\frac{s}{s-\rho}+\frac{\sigma^2}{2}\frac{\xi^{'}(0)}{s-\rho}+\frac{c}{s-\rho}-\frac{\lambda r\hat{w_d}(s)}{s-\rho}.
\end{equation}
Dividing equation (\ref{glap}) by $\frac{\sigma^2}{2}(s+\rho)+c$ and rearranging terms lead to:
\begin{equation}\label{glapl}
    \hat{\xi}(s)=\frac{2\lambda r/\sigma^2}{s+\rho+2c/\sigma^2}\hat{\xi}(s)T_sT_\rho f(0)+\frac{\rho+2c/\sigma^2+\xi^{'}(0)}{(s-\rho)(s+\rho+2c/\sigma^2)}+\frac{1}{s+\rho+2c/\sigma^2}-\frac{2\lambda r/\sigma^2\hat{w_d}(s)}{(s-\rho)(s+\rho+2c/\sigma^2)}.
    \end{equation}
Inverting the Laplace transforms in (\ref{glapl}) gives:
\begin{eqnarray*}\label{gx}
                \xi(x) &=& \frac{2\lambda r}{\sigma^2}\int_0^{x}\xi(x-y)\beta*T_\rho f(y)dy+[\rho+2c/\sigma^2+\xi^{'}(0)]f*\beta(x)+\beta(x)-\frac{2\lambda r}{\sigma^2}f*\beta*w_d(x)\nonumber\\
                       &=& \sum_{n=0}^{\infty}\left(\frac{2\lambda r}{\sigma^2}\right)^{n}(\beta*T_{\rho}p)^{*n}*\varphi(x).
                        \end{eqnarray*}
 Thus
\begin{eqnarray}\label{gx}
                h^d(x) &=& h^d(0)\xi(x)
                       = h^d(0)\left[ \sum_{n=0}^{\infty}\left(\frac{2\lambda r}{\sigma^2}\right)^{n}(\beta*T_{\rho}p)^{*n}*\varphi(x)\right].
 \end{eqnarray}
From the boundary condition (\ref{bhd})
we get that:
\begin{equation*}
h^d(0)=\frac{1}{\sum_{n=0}^{\infty}(\frac{2\lambda r}{\sigma^2})^{n}(\beta*T_{\rho}p)^{*n}*\varphi(a)},
\end{equation*}
and by (\ref{gx}) we derive (\ref{hdxr}).
The cases $\sigma=0$ and $d=0$ could be analyzed analogously.
This completes the proof.
\vspace{3mm}\hfill $\Box$

\section{Optimality}
In this section we will give sufficient conditions for the barrier strategy $\pi^{a^*}$ to be optimal.
We will now give the main results of this paper which give sufficient conditions for optimality of the barrier strategy $\pi^{a^*}$.

Let $\Gamma$ be the generator of $X_t$ acting in sufficiently smooth functions $g$, defined by
\begin{equation*}
\Gamma g(x)=cg^{'}(x)-\lambda g(x)+\lambda r\int_{(0,\infty)}g(x-y)f(y)dy.
\end{equation*}

We will start from classical Verification Lemma.
\begin{lem}\label{verlem}
Suppose $\pi$ is an admissible dividend strategy such that $v^{\pi}$ is sufficiently smooth and for all $x>0$:
\begin{equation}\label{HJB}
\max\{(\Gamma-q)v^{\pi}(x), 1-(v^{\pi})^{'}(x)\}\leq0.
\end{equation}
Then $v^{\pi}(x)=v^*(x)$ for all $x\in\mathbb{R}$.
\end{lem}
{\bf Remark.} The inequality (\ref{HJB}) it is called Hamiltonian-Jacobi-Bellman system
and is classical in all stochastic optimization problems.

\textbf{Proof.}
We will follow classical arguments. Suppose that $g$ is $C^2(\mathbb{R})$ and that
$$\max\{(\Gamma-q)g(x), 1-g^{'}(x)\}\leq0.$$
We will prove that
\begin{equation}\label{HJBb}
g(x)\geq\sup_{\pi\in\Pi} v^\pi(x),\qquad x\in \mathbb{R}.
\end{equation}
Having strategy $\pi^*$ for which $g(x)=v^{\pi^*,d}(x)$ will complete the proof.
To prove (\ref{HJBb}) we will consider
Markov process $(t,N_t,X_t, \varsigma^X_t)$, with $\varsigma^Z_t
=t-\sup\{s\leq t: Z_{t}\geq 0\}$ for some process $Z$.
By Sato \cite[Ch. 6, Thm. 31.5]{Sato}
function $\varpi(t,k,x, z):=r^ke^{-qt}g(x)\mathbf{1}_{\{z\leq r\}}$ is in a domain of extended generator $\mathcal{A}$
of this four-dimensional process.
In fact using similar arguments like in deriving euation (\ref{hidea})
from the definition of the infinitesimal generator one can prove that
$$(\mathcal{A}-q\mathbb{I})\varpi(t,k,x, z)=r^k(\Gamma -q)g(x).$$
Recall that $U_t^\pi=X_t-L_t^\pi$ and note that regulations $L_t^\pi$ do not modify jump process of $X$.
Note that finite number of discontinuities
in $g$ and hence also single discontinuity in $\mathbf{1}_{\{z\leq r\}}$ are allowed here.
Hence we are also allowed
to apply It\^{o}'s lemma and letting below $U=U^\pi$ and $L=L^\pi$ we derive:
\begin{eqnarray}
\mathrm{e}^{-qt}r^{N_t}g(U_t)\mathbf{1}_{\{\varsigma^U_t\leq r\}} - g(U_0) &=& r^{N_t}J_{g}(t) - r^{N_t}\int_0^t\mathrm{e}^{-qs}g^\prime(U_{s^-})d L^c_s\nonumber\\
&&+ \int_0^t\mathrm{e}^{-qs}r^{N_s}(\Gamma g - qg)(U_{s-})d s + M_t,
\label{eq:ITO}
\end{eqnarray}
where $M_t$ is a local martingale with $M_0=0$, $L^c$ is the
pathwise continuous parts of $L$
and
\begin{equation}\label{eq:sumJ}
J_g(t) = \sum_{s\leq t}\mathrm{e}^{-qs}\left[ g(A_s + B_s) - g(A_s)\right]\mathbf{1}_{\{B_s\neq 0\}},
\end{equation}
where $A_s = U_{s-}+\Delta X_s$ with $\Delta X_s=X_s-X_{s-}$ and $B_s = -\Delta L_s$ denotes the jump
of $-L$ at time $s$. Let $T_n$ be a localizing sequence of $M$.
Applying the Optional Stopping Theorem to the stopping times
$T_k^\prime=T_k\wedge \tau^{\pi,d}$ and using the Fatou's Theorem  we
derive:
\begin{eqnarray*}
g(x)&\geq& \mathbb{E}_xr^{N_{T^\prime_n}}\left[\mathrm{e}^{-qT_n^\prime}g(U_{T^\prime_n})\mathbf{1}_{\{\varsigma^U_{T^\prime_n}\leq r\}}-
J_g(T^\prime_n)\right]+
\mathbb{E}_xr^{N_{T^\prime_n}}\int_0^{T^\prime_n}\mathrm{e}^{-qs}g^\prime(U_{s^-})d L^c_s
\\\\
&&- \mathbb{E}_x\int_0^{T^\prime_n}r^{N_s}\mathrm{e}^{-qs}(\Gamma g - qg)(U_{s-})d s.
\end{eqnarray*}
Invoking the variational inequalities $g^\prime(x)\geq 1$ (hence $g(A_s + B_s) - g(A_s)\leq - \Delta L_s$ if
$A_s>0$) and $(\Gamma-q)g(x)\leq 0$ we have:
\begin{eqnarray*}g(x)&\geq& \mathbb{E}_x\mathrm{e}^{-qT_n^\prime}r^{N_{T^\prime_n}}g(U_{T^\prime_n})\mathbf{1}_{\{\varsigma^U_{T^\prime_n}\leq d\}}
+\mathbb{E}_xr^{N_{T^\prime_n}}\int_0^{T^\prime_n}\mathrm{e}^{-qs}d L_s\\&\geq&
\mathbb{E}_x\left[r^{N_{\tau^{\pi,d}}}\mathrm{e}^{-q\tau^{\pi,d}}g(U_{\tau^{\pi,d}}); \tau^{\pi,d}\leq
T_n^\prime\right] +\mathbb{E}_xr^{N_{\tau^{\pi,d}}}\left[\int_0^{\tau^{\pi,d}}\mathrm{e}^{-qs}d L_s;
\tau^{\pi,d}\leq T_n^\prime\right].
\end{eqnarray*}
Letting $n\to\infty$ in the conjunction with the monotone convergence
theorem and using the fact that $\mathbf{1}_{\{\varsigma^U_{\tau^{\pi,d}}\leq d\}}=0$ complete the proof.
\vspace{3mm}\hfill $\Box$

We will now focus on the optimality of the barrier strategy $\pi^{a^*}$.
\begin{theo}\label{fin}
If
\begin{equation}\label{condition1}
(\Gamma -q)h^d(x)\leq  0\qquad \text{for $x\geq a^*$}
\end{equation}
then the barrier strategy $\pi^{a^*}$ is optimal.
\end{theo}
\textbf{Proof.}
Note that from (\ref{hide}) it follows that
\begin{equation}\label{zeromart}
(\Gamma -q)h^d(x)=  0\quad \text{for $x\leq a^*$}.
\end{equation}
From the choice of $a^*$ and Theorem \ref{jieguoteo}
it follows that $(v^{a^*,d})'(x)\leq 1$.
This completes the proof in view of the Verification Lemma \ref{verlem}.
\vspace{3mm}\hfill $\Box$

Moreover, we can give other necessary condition for the barrier strategy to be optimal.
First note that the function $h^d(x)$ could be formally defined via $g^d(x)$ (and hence via general upper level $b$)
by (\ref{twosidedexit}) for any $x\in\mathbb{R}$.
Moreover we can exchange $h^d$ into $g^d$ in the value function $v^{d,a^*}$ given in
Theorem \ref{jieguoteo}.

\begin{coro}\label{wniosek}
Suppose that
$$(g^{d})'(a)\leq
(g^{d})'(b),\qquad\mbox{for all $a^*\leq a\leq b$.}
$$
Then the barrier strategy $\pi_{a^*}$ at $a^*$ is an optimal strategy.
\end{coro}
\textbf{Proof.}
Using the Theorem \ref{fin} and (\ref{zeromart}) the proof of this result is the same as the proof of \cite[Theorem 2]{loeffen2008optimality}.\
\vspace{3mm}\hfill $\Box$

\begin{theo}\label{monotone}
Suppose that the claims $C_i$ ($i=1,2\ldots$) have density with $f'(y)$ being monotonne decreasing.
Then $v^{a^*,d}(x)=v^*(x)$ and $\pi_{a^*}$ is an optimal strategy.
\end{theo}

\textbf{Proof.}
We start from the basic observation that
process $\alpha X_t+\log rN_t$ is a L\'evy process (in fact is a compund Poisson process with drift perturbed by
independent Brownian motion). Hence from L\'evy-Khinchine theorem
there exists function $\psi_r(\alpha)$ (which is discounted by the number of jumps Laplace exponent of $X$)
such that
$$\log \mathbb{E}r^{N_t}e^{\alpha X_t}=\log \mathbb{E}e^{\alpha X_t+N_t\log r}=\psi_r(\alpha)t.$$
Now we can define new probability measure:
\begin{equation}
\label{eq:changeofmeasure} \left.\frac{d \mathbb{P}^\alpha}{d\mathbb{P}}\right|_{\mathcal{F}_t} =
r^{N_t}\exp\left(\alpha X_t - \psi_r(\alpha)t\right),
\end{equation}
where ${\mathcal{F}_t}$ is a right-continuous natural filtration  of $X$.
Under the measure $\mathbb{P}^\alpha$ the process $X$ is still in the same form with the new parameters
(see Palmowski and Rolski \cite{PR}).
It easy to prove (like in the case of spectrally negative L\'evy processes) that
for fixed $r$ function $\alpha\to\psi_r(\alpha)$ is convex and there exists $\Phi_r(q)$ such that
$$\psi_r(\Phi_r(q))=q.$$
Using optional stopping theorem with the observation that $X(\tau_b^+)=b$ and the fact that $g^d(x)$ could be defined
for any $b>0$ (hence also $b=+\infty$) we derive new representation:
\begin{equation}\label{gdanew}
g^d(x)=e^{\Phi_r(q)x}\mathbb{P}^{\Phi_r(q)}_x(\tau^d=\infty)
\end{equation}
Now the proof follows the same idea like the proof of Corollary
5.2 of Czarna and Palmowski \cite{czarna2010dividend}.
\vspace{3mm}\hfill $\Box$
\begin{re}\rm
From the equation (\ref{twosidedexit}) and
definition of function $h^d(x)$ given in (\ref{hddefinition}) we also proved another interesting for itself identity:
$$\mathbb{E}_x[r^{N_{\tau^+_a}}e^{-q\tau^+_a}, \tau^+_a<\tau^d]
=e^{\Phi_r(q)(x-a)}\frac{\mathbb{P}^{\Phi_r(q)}_x(\tau^d=\infty)}{\mathbb{P}^{\Phi_r(q)}_a(\tau^d=\infty)}$$
giving the solution of Parisian-type two-sided exit problem discounted by the number of claims;
compare with equation (18) of Czarna and Palmowski \cite{czarna2010dividend}.
\end{re}

\section{Exponential claims}

Now consider the case when claim sizes are exponentially distributed with parameter $\mu$, that is
$f(x)=\mu e^{-\mu x}$ for $x>0$ and $\sigma=0$.

\textbf{Intoductionary expressions.}

Utilizing the fact that $f$ is an exponential density $T_{\rho}f(x)$ becomes
\begin{eqnarray*}
T_{\rho}f(x)&=&\int_{x}^{\infty}e^{-\rho(y-x)}\mu e^{-\mu y}dy\nonumber\\
&=&\frac{\mu}{\rho+\mu}e^{-\mu x}.
\end{eqnarray*}
The $n$-fold convolution of $T_{\rho}f$ becomes
\begin{eqnarray}\label{expp}
(T_{\rho}f)^{*n}(x)&=&\left(\frac{\mu}{\rho+\mu}\right)^n\frac{x^{n-1}}{(n-1)!}e^{-\mu x}, \ \ \ \ n>1
\end{eqnarray}
 and $w(x)$ becomes
 \begin{equation}\label{expw}
 w(x)=e^{-\mu x}u(d),
  \end{equation}
 where
 \[u(d)=\sum_{k=0}^{\infty}\frac{r^k\lambda^k(\mu c)^{k+1}}{k!(k+1)!}\int_0^dt^{2k}e^{-(\lambda+q+\mu c)t}dt.\]

  Substituting (\ref{expp}) and (\ref{expw}) into (\ref{hxr}), (\ref{hdxr}), respectively, for $0\leq x\leq a$
  gives:
  \begin{equation}\label{ehdx}
h^{d}(x)=\left\{
  \begin{array}{lll}
&\frac{\vartheta(x)}{\vartheta(a)},&\ d=0,\\
&\frac{\varrho(x)}{\varrho(a)},& \ d>0,
  \end{array}
\right.
\end{equation}
 where
\begin{eqnarray*}
\vartheta(x)&=&\sum_{n=1}^{\infty}\frac{(\lambda r \mu)^n}{(n-1)!c^n(\rho+\mu)^{n}}e^{\rho x}\int_0^xy^{n-1}e^{-(\rho+\mu)y}dy+e^{\rho x},\nonumber\\
\varrho(x)&=&\sum_{n=1}^{\infty}\frac{(\lambda r \mu)^n}{(n-1)!c^n(\rho+\mu)^{n}}\left[(1-\frac{\lambda r u(d)}{c(\rho+\mu)})e^{\rho x}\int_0^xy^{n-1}e^{-(\rho+\mu)y}dy+\frac{\lambda r, u(d)}{cn(\rho+\mu)}x^ne^{-\mu x}\right]\nonumber\\
&&+e^{\rho x}-\frac{\lambda r u(d)}{c(\rho+\mu)}(e^{\rho x}-e^{-\mu x}).
\end{eqnarray*}
Moreover,
 \begin{equation}\label{mde}
  \vartheta^{'}(x)=\sum_{n=1}^{\infty}\frac{(\lambda r \mu)^n}{(n-1)!c^n(\rho+\mu)^{n}}\left[\rho e^{\rho x}\int_0^xy^{n-1}e^{-(\rho+\mu)y}dy+x^{n-1}e^{-\mu x}\right]+\rho e^{\rho x}
  \end{equation}
  and
\begin{eqnarray}\label{nde}
\varrho^{'}(x)&=&\sum_{n=1}^{\infty}\frac{(\lambda r \mu)^n}{(n-1)!c^n(\rho+\mu)^{n}}\left[(1-\frac{\lambda r u(d)}{c(\rho+\mu)})\rho e^{\rho x}\int_0^xy^{n-1}e^{-(\rho+\mu)y}dy+x^{n-1}e^{-\mu x}-\frac{\lambda\mu r u(d)}{cn(\rho+\mu)}x^ne^{-\mu x}\right]\nonumber\\
&&+\rho e^{\rho x}-\frac{\lambda r u(d)}{c(\rho+\mu)}(\rho e^{\rho x}+\mu e^{-\mu x}).
\end{eqnarray}
Further, taking derivative for $\vartheta^{'}(x)$ and $\varrho^{'}(x)$, we can get the expressions of $\vartheta^{''}(x)$ and $\varrho^{''}(x)$.
Then solving equations $\vartheta^{''}(x)=0$ and $\varrho^{''}(x)=0$ we will derive the optimal barrier levels $a^*$ and $b^*$ for $d=0$ and $d>0$, respectively.

Recall that for $d=0$,
\begin{equation}\label{jieguo5}
v^{a^*,d}(x)=\left\{
  \begin{array}{lll}
&\frac{\vartheta(x)}{\vartheta^{'}(a^*)}&\text{if}\ 0\leq x\leq a^*,\\
&x-a^*+\frac{\vartheta(a^*)}{\vartheta^{'}(a^*)}&\text{if} \ x> a^*
  \end{array}
\right.
\end{equation}
and for $d>0$,
\begin{equation}\label{jieguo6}
v^{b^*,d}(x)=\left\{
  \begin{array}{lll}
&\frac{\varrho(x)}{\varrho^{'}(b^*)}&\text{if}\ x\leq b^*,\\
&x-b^*+\frac{\varrho(b^*)}{\varrho^{'}(b^*)}&\text{if} \ x> b^*,
  \end{array}
\right.
\end{equation}
where $a^*$ is the optimal barrier when $d=0$ and $b^*$ is the optimal barrier when $d>0$.

\textbf{Numerical analysis.}

Using above expressions we were able to find function $h^d(x)$ numerically and hence
the value function. We also identified the optimal barriers.
This shows that produced algorithms can be used in daily practice.

Let $\lambda=10, \mu=1, c=15, q=0.1,r=0.8$. We will consider two cases: $d=0$ and $d=2$.
The positive root of equation (\ref{lundberg}) is $\rho=0.24493$.
For $d=0$ solving equation $\vartheta^{''}(x)=0$ produces $a^*\thickapprox0.7693$ and for
$d=2$ solving equation $\varrho^{''}(x)=0$ gives $b^*\thickapprox0.52202$.
In the following Figures 3-4 we show how $h^d$ looks like.
With the help of Mathematica we also plotted $(\Gamma-q)v^{a^*,d}$  and $(\Gamma-q)v^{b^*,d}$
(see Figures 3-4).
In all cases we have $(\Gamma-q)v^{a^*,d}\leq 0$, hence by Theorem \ref{fin} the barriers $a^*$ and $b^*$ produce optimal strategies
in these both cases (it follows also straightforward from Theorem \ref{monotone} and the choice of the exponential claim size).

\begin{figure}[htbf]
 \centering
 \begin{minipage}[t]{0.45\textwidth}
 \centering
 \includegraphics[scale=0.5]{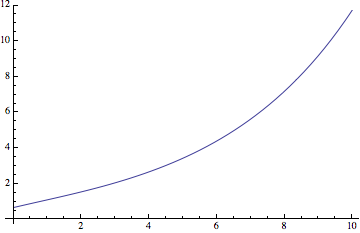}
 \caption*{ $h^d(x)$}
  \end{minipage}
  \begin{minipage}[t]{0.45\textwidth}
   \centering
   \includegraphics[scale=0.5]{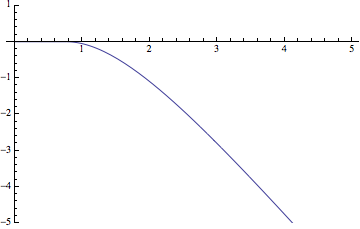}
 \caption*{ $(\Gamma-q)v^{a^*,d}$}
  \end{minipage}
   \caption*{Fig.3.   The graphs of $h^d(x)$ and $(\Gamma-q)v^{a^*,d}$ for $\lambda=10, \mu=1, c=15, q=0.1,r=0.8$ and $d=0.$ }
  \end{figure}

\begin{figure}[htbf]
 \centering
 \begin{minipage}[t]{0.45\textwidth}
 \centering
 \includegraphics[scale=0.5]{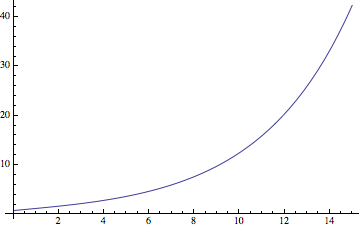}
 \caption*{ $h^d(x)$}
  \end{minipage}
  \begin{minipage}[t]{0.45\textwidth}
   \centering
   \includegraphics[scale=0.5]{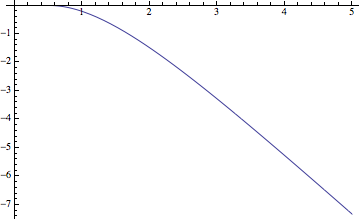}
 \caption*{ $(\Gamma-q)v^{b^*,d}$}
  \end{minipage}
   \caption*{Fig.4.  The graphs of $h^d(x)$ and $(\Gamma-q)v^{b^*,d}$ for $\lambda=10, \mu=1, c=15, q=0.1,r=0.8$ and $d=2.$ }
  \end{figure}

\noindent{\bf{Acknowledgements}}\\
This research is support by the FP7 Grant PIRSES-GA-2012-318984.
Zbigniew Palmowski is supported by the Ministry of Science and Higher Education of Poland under the grant
2013/09/B/HS4/01496.

\bibliographystyle{plain}

\bibliography{refbib}
\end{document}